\documentclass[10pt]{article}
\usepackage{amssymb}
\usepackage{latexsym,bm}
\usepackage{graphicx}
\usepackage{mathrsfs,amscd,amssymb,amsthm,amsmath,bm,graphicx,psfrag,subfigure,url,xcolor}
\usepackage{amsmath}
\usepackage{mathrsfs}

\date{}
\begin{document}
\title{ Extremal octagonal chains with respect to the Kirchhoff index \footnote{E-mail addresses:
{\tt mathqima@163.com}(Q. Ma).}}
\author{\hskip -10mm  Qi Ma\thanks{Corresponding author.}\\
{\hskip -10mm \small Sanda University, Shanghai 200241, China}}\maketitle
\begin{abstract}
Let $G$ be a connected graph. The resistance distance between any two vertices of $G$ is equal to the effective resistance between them in the corresponding electrical network constructed from $G$ by replacing each edge with a unit resistor. The Kirchhoff index is defined as the sum of resistance distances between all pairs of the vertices. These indices have been computed for many interesting graphs, such as linear polyomino chain, linear/M\"{o}bius/cylinder hexagonal chain, and linear/M\"{o}bius/cylinder octagonal chain. In this paper, we characterized the maximum and minimum octagonal chains with respect to the Kirchhoff index.
\end{abstract}


\vskip 8mm

\section{ Introduction}
In this paper, we consider simple connected graph $G=(V(G),E(G))$ with vertex set $V(G)$ and edge set $E(G)$. We call $|V(G)|$ the order of $G$ and $| E(G )|$ the size of $G.$  For two vertices $u$ and $v$, we use the symbol $u\leftrightarrow v$ to mean that $u$ and $v$ are adjacent and use $u\nleftrightarrow v$ to mean that $u$ and $v$ are non-adjacent. Let $d_G(u,v)$ be the distance between vertices $u$ and $v$ in $G$, which represents the length of a shortest path connecting vertex $u$ and $v$ in $G$. We follow the book \cite{7} for terminology and notations.

Molecules can be modeled by graphs with vertices for atoms and edges for atomic bonds. The topological indices of a molecular graph can provide some information on the chemical properties of the corresponding molecule. It plays essential roles in the study of QSAR/QSPR in chemistry. The {\it Wiener index} of $G$, introduced in \cite{8}, is defined as
$$
 W(G)=\sum_{\{u,v\}\subseteq V(G)}d_G(u,v).
$$
Wiener \cite{8} used it to study the boiling point of paraffin. Many chemical properties of molecules are related to the Wiener index. Based on the electronic network theory, Klein and Randi\'{c} \cite{9} proposed the concept of {\it resistance distance} in 1993. The resistance distance $r_G(u,v)$ between vertices $u$ and $v$ of a connected ( molecular) graph $G$ is computed as the effective resistance between vertices $u$ and $v$ in the corresponding electrical network constructed from $G$ by replacing each edge of $G$ with an unit resistor. Without causing ambiguity, we replace $r_G(u,v)$ with $r(u,v)$. This novel parameter is in fact intrinsic to the graph and has some nice interpretations and applications in chemistry. Similar to the Wiener index,  Klein and Randi\'{c} \cite{7} defined the {\it Kirchhoff index} $K\!f(G)$ of $G$ as the sum of the pairwise resistance distances between vertices, i.e.
$$
 K\!f(G)=\sum_{\{u,v\}\subseteq V(G)}r_G(u,v).
$$

In this paper, we focus only on octagonal chains. Let $Q_{n}$ be a ladder graph (a linear quadrilateral chain) with $n$ squares. Denote by $S_i,i=1,\ldots,n$ the $i$-th square of $Q_{n}$. Note that a octagonal chain $O_n$ with $n$ octagons can be obtained from $Q_{n}$ by adding four vertices to each of the $1$-st, $2$-rd,\ldots, $n$-th square. Obviously, We have five ways to add these four new vertices to $S_{i},i=1,2,\ldots,n.$ That is, we can add $0$ (resp. $1, 2, 3$ or $4$) vertices to the top edge of $S_{i}$ and the remaining vertices to the bottom edge of $S_{i}.$ For convenience, we always suppose that we add four vertices to the bottom edge of $S_1$ and $S_{n}$. For each of the remaining $(n-2)$-octagons, we give a number $w_i=0$ (resp. $1, 2, 3$ or $4$) to the octagon if the octagon is obtained by adding $0$ (resp. $1, 2, 3$ or $4$) vertex to the top edge. In this viewpoint, we are able to represent a octagonal chain with $n$ octagons by a $(n-2)$-vector $w=(w_1, w_2,\ldots ,w_{n-2})$ such that $w_i\in \{0,1,2,3,4\}$. In the following, we always denote a octagonal chain with $n$ octagons by $O(w)$ such that $w$ is a $(n-2)$-tuple of $0, 1, 2, 3$ or $4$. A {\it kink} in a octagonal chain is a octagon whose $w_i=0$ or $4.$ A octagonal chain $G(w)$ with $w_i=0$ or $4$ is called a {\it ``all-kink" chain}, where $1\leq i\leq n-2$. There are two special octagonal chains. We call $O(\underbrace{2,2,...,2}_{n-2})$ a {\it linear octagonal chain}, and denote it by $L_n.$ The octagonal chain $O(\underbrace{0,0,...,0}_{n-2})$ is called a {\it helicene octagonal chain}, which is denoted by $D_n$. $L_5$ and $D_5$ are illustrated in Fig. 1.
\begin{figure}[!ht]
  \centering
 \includegraphics[width=100mm]{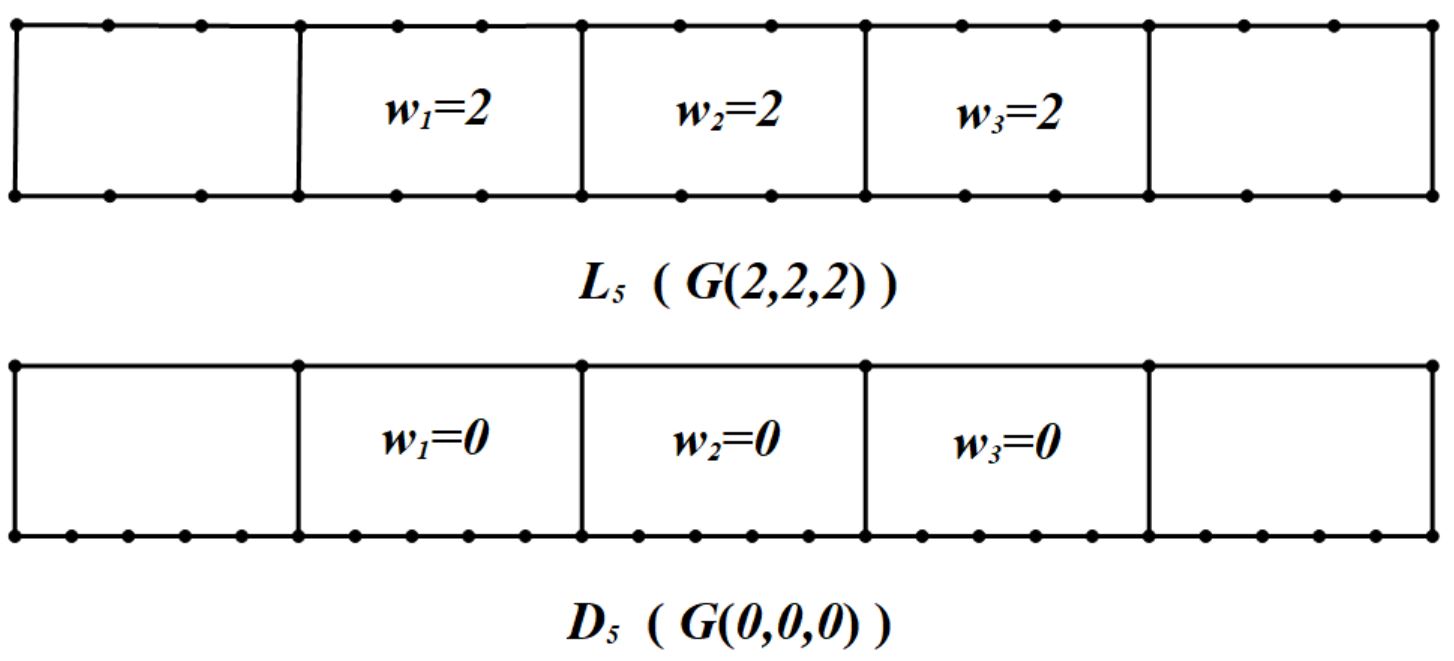}
 \caption{Some phenylene chains with $5$ octagons.}
\end{figure}

In this paper, we consider Kirchhoff indices of an interesting class of molecular graphs, octagonal chains, which are molecular graphs of an important class of unbranched conjugated hydrocarbons. In 2020, Zhao, Liu and Hayat \cite{10} determined the Resistance distance-based graph invariants and the number of spanning trees of linear crossed octagonal graphs.  Liu, Zhang, Wang and  Lin \cite{11} determined the Kirchhoff index and spanning trees of M\"{o}bius/cylinder octagonal chain. Zhu \cite{12} studied the Kirchhoff index, degree-Kirchhoff index and spanning trees of linear octagonal chains. Geng ,Li and Wei identified the extremal octagonal chains with n octagons having the maximum and minimum spectral radii.  More interesting results about these topics can be found in \cite{13,14}. In these paper, we determine the unique chain with maximum Kirchhoff index and the unique chain with minimum Kirchhoff index. Our main results as follows:

{\bf Theorem 1.} {\it Among all octagonal chains, the linear octagonal chain is the unique chain with maximum Kirchhoff index. }

{\bf Theorem 2.} {\it Among all octagonal chains, the helicene octagonal chain is the unique chain with minimum Kirchhoff index. }

\section{Preliminaries}
To prove our main result, we will need the following lemmas and definitions.

{\bf Definition 1.} (Series Transformation) Let $x, y$ and $z$ be nodes in a graph where $y$ is adjacent to only $x$ and $z$. Moreover, let $R_1$ equal the resistance between $x$ and $y$ and $R_2$ equal the resistance between node $y$ and $z$. A series transformation transforms this graph by deleting $y$ and setting the resistance between $x$ and $z$ equal to $R_1 + R_2$.

{\bf Definition 2.} (Parallel Transformation) Let $x$ and $y$ be nodes in a multi-edged graph where $e_1$ and $e_2$ are two edges between $x$ and $y$ with resistances $R_1$ and $R_2$, respectively. A parallel transformation transforms the graph by deleting edges $e_1$ and $e_2$ and adding a new edge between $x$ and $y$ with edge resistance $r=(\frac{1}{R_1}+\frac{1}{R_2})^{-1}$

A $\Delta$-$Y$ transformation is a mathematical technique to convert resistors in a triangle formation to an equivalent system of three resistors in a $``Y"$ format as illustrated in Fig. 2. We formalize this transformation below.
\begin{figure}[!ht]
\centering
  \includegraphics[width=80mm]{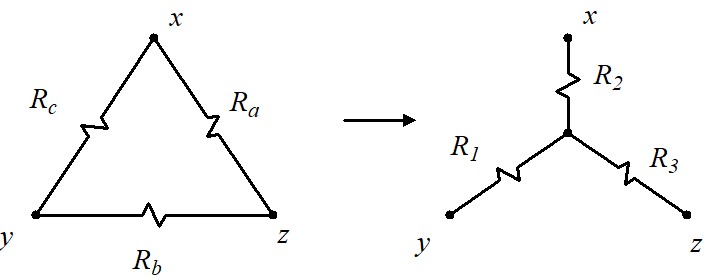}
  \caption{$\Delta$ and $Y$ circuits with vertices labeled as in Definition 3.}
\end{figure}

{\bf Definition 3.} ($\Delta$-$Y$ Transformation) Let $x,y,z$ be nodes and $R_a,R_b$ and $R_c$ be given resistances as shown in Fig. 2. The transformed circuit in the $``Y"$ format as shown in Fig. 2 has the following resistances:
\begin{eqnarray*}
  R_1=\frac{R_bR_c}{R_a+R_b+R_c}, \ \  R_2=\frac{R_aR_c}{R_a+R_b+R_c}, \ \
  R_3=\frac{R_aR_b}{R_a+R_b+R_c}.
\end{eqnarray*}

{\bf Lemma 3.} (Stevenson \cite{6}) {\it Series transformations, parallel transformations, and $\Delta-Y$ transformations yield equivalent circuits.}

The concept of $S,T$-isomers was introduced by Polansky and Zander \cite{5} in 1982. From then on, a lots of research \cite{1},\cite{2}, \cite{3}, \cite{4} has been devoted to the study of topological properties of $S,T$-isomers. Let $S$ be a graph with edge cut $[A,B]=\{ux,\,vy\}$ as shown in Fig. 3. Note that $u$ and $u$ are two distinct vertices of $A$ and $x$ and $y$ are two distinct vertices of $B$. The graph $T$ is obtained from $S$ by deleting edges $\{ux,\,vy\}$ and adding edges $\{uy,\,vx\};$ see Figure 3.
\begin{figure}[!ht]
\centering
  \includegraphics[width=60mm]{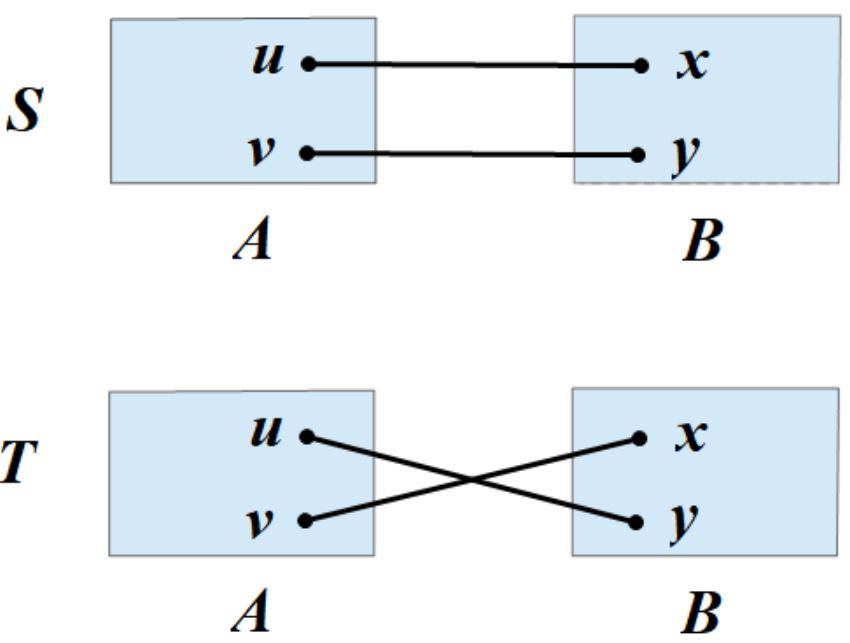}
  \caption{The structure of the graphs $S$ and $T$ and the labeling of their vertices.}
\end{figure}

{\bf Definition 4.} In the following, we will use $r_G(x)$ to denote the sum of resistance distances between $x$ and all the other vertices of $G.$ More precisely,
$$
r_G(x)=\sum_{y\in V(G)\setminus\{x\}}r_G(x,y).
$$
Yang and Klein \cite{2} obtained the comparison theorem on the Kirchhoff index of $S,T$-isomers, which plays an essential role in the characterization of extremal octagonal chains.

{\bf Lemma 4.} (\cite{2}) {\it Let $S,T,A,B,u,v,x,y$ be defined as in Fig. 3. Then}
$$
K\!f(S)-K\!f(T)=\frac{[r_A(u)-r_A(v)][r_B(y)-r_B(x)]}{r_A(u,v)+r_B(x,y)+2}.
$$

\section{Proof of the Main Result}
In order to determine extremal octagonal chains with respect to Kirchhoff index, we first give some comparison results on Kirchhoff index in octagonal chains. Recalled that we use $O(w)$ to denote an octagonal chain with $n$ octagons and $w=(w_1,w_2,\ldots,w_{n-2}).$
 \begin{figure}[!ht]
\centering
  \includegraphics[width=100mm]{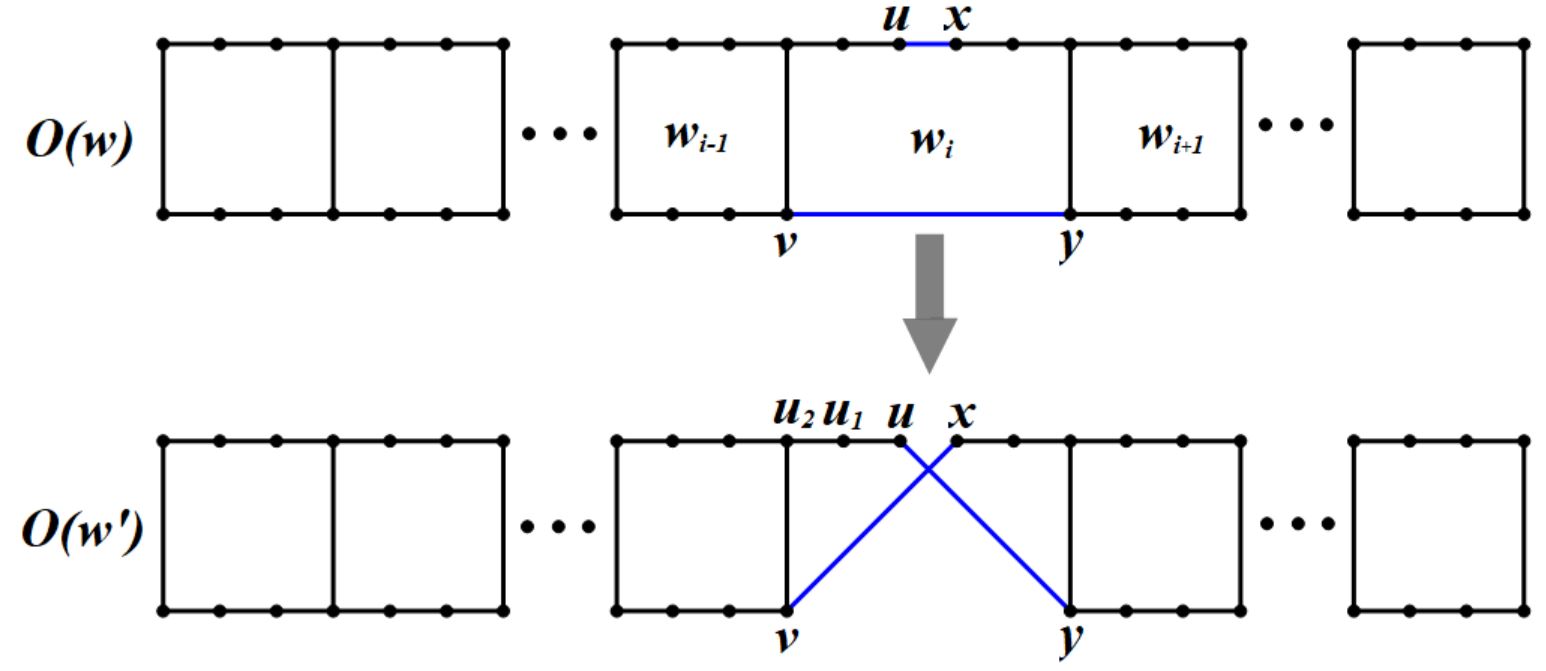}
  \caption{$O(w)$ and $O(w')$ in the proof of Lemma 5.}
\end{figure}

{\bf Lemma 5.} {\it If $O(w)$ is a octagonal chain with some integer $i\in\{1,2,\ldots,n-2\}$ such that $w_i=0$ or $4,$ then let $w'=\{w_1,\ldots,w_{i-1},2,4-w_{i+1},\ldots,4-w_{n-2}\},$
and we have $K\!f(O(w))<K\!f(O(w')).$}

{\bf Proof.} We first prove the case $w_i=4$. The proof of $w_i=0$ is similar to the proof of $w_i=4$, and we omit the process. First choosing four vertices $u,v,x$ and $y$ in the $(i+1)$-th octagon of $O(w)$ as shown in Fig. 4. Note that edges set $\{ux\,,vy\}$ is a edge cut of $O(w).$ Deleting edges $\{ux\,,vy\}$ and adding edges $\{uy\,,vx\},$ we can get graph $O(w').$ By the construction of $O(w')$, we deduce that $O(w)$ and $O(w')$ are pairs of $S,T$-isomers. Let $A_1$ be the component of $O(w)-\{ux\,,vy\}$ such that $\{u,v\}\in V(A_1)$ and let $B_1$ be the component of $O(w)-\{ux\,,vy\}$ such that $\{x,y\}\in V(B_1).$ Denote by $u_1$ the unique neighbor of $u$ in $A_1.$ Note that $u_1$ has degree two. We denote by $u_2$ the other neighbor of $u_1;$ See Fig. 4. Note that if $k\in V(A_1)\backslash \{u_1,u\},$ we have $r_{A_1}(u,k)=r_{A_1}(u_2,k)+2.$

By Lemma 4, We have
$$
K\!f(O(w))-K\!f(O(w'))=\frac{[r_{A_1}(u)-r_{A_1}(v)][r_{B_1}(y)-r_{B_1}(x)]}{r_{A_1}(u,v)+r_{B_1}(x,y)+2}.
$$
Now, we are going to compare $r_{A_1}(u)$ and $r_{A_1}(v)$. We distinguish two cases. If $k\in V(A_1)\backslash \{u,u_1,u_2,v\},$ by the triangular inequality, we have
$$
r_{A_1}(v,k)\leq r_{A_1}(v,u_2)+r_{A_1}(u_2,k)< 1+r_{A_1}(u_2,k)< r_{A_1}(u,k),
$$
where the second inequality follows from the fact that $r_{A_1}(v,u_2)<d_{A_1}(v,u_2)=1.$  If $k\in \{u_1,u_2\},$ we have
$$
r_{A_1}(v,u_1)+r_{A_1}(v,u_2)=1+2r_{A_1}(v,u_2)< 3 =r_{A_1}(u,u_1)+r_{A_1}(u,u_2).
$$
Since
$$
r_{A_1}(u)-r_{A_1}(v)=\sum_{k\in V(A_1)\backslash\{u,v\}}r_{A_1}(u,k)-\sum_{k\in V(A_1)\backslash\{u,v\}}r_{A_1}(v,k),
$$
we deduce that $r_{A_1}(u)-r_{A_1}(v)>0.$ By the same argument, we can obtain $r_{B_1}(y)-r_{B_1}(x)<0.$ Thus
$$
K\!f(O(w))-K\!f(O(w'))=\frac{[r_{A_1}(u)-r_{A_1}(v)][r_B(y)-r_B(x)]}{r_{A_1}(u,v)+r_{B_1}(x,y)+2}<0.
$$
This completes the proof. \hfill $\Box$

\begin{figure}[!ht]
\centering
  \includegraphics[width=100mm]{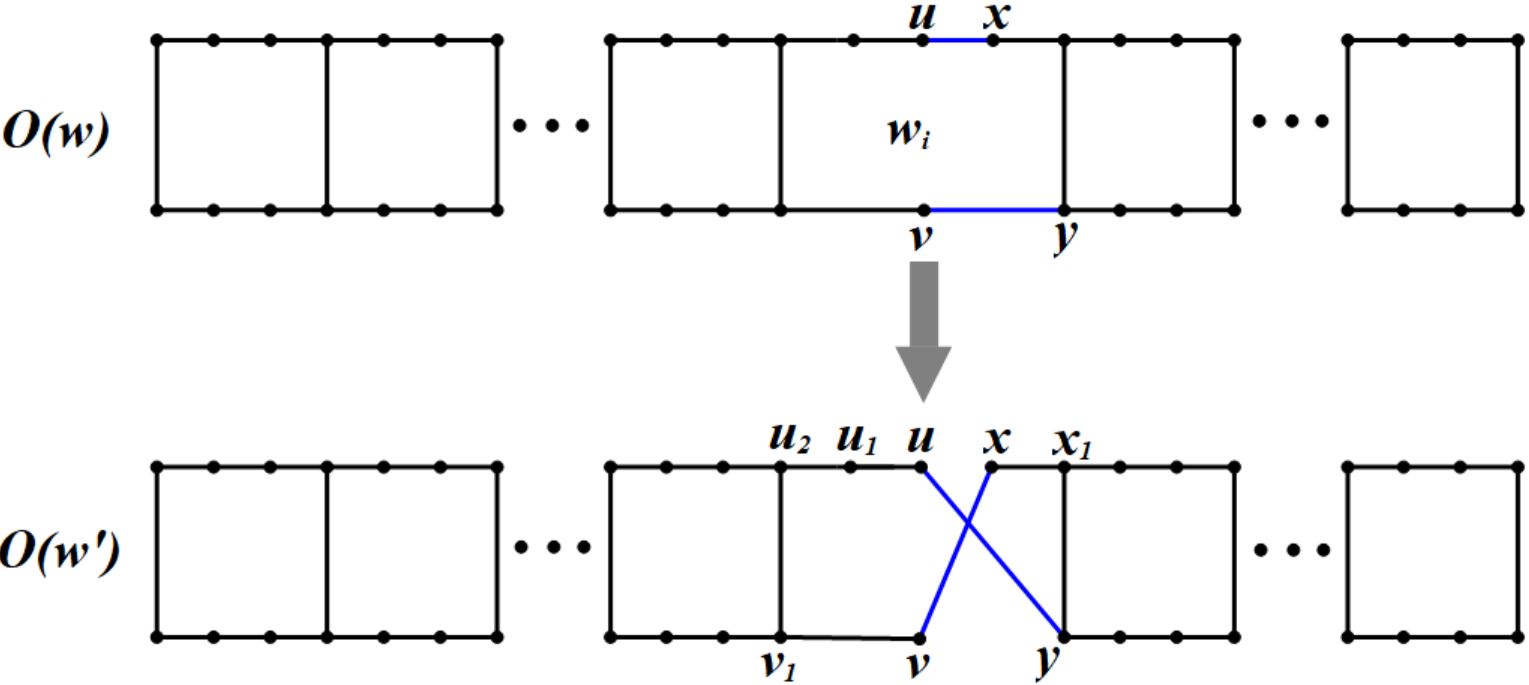}
  \caption{$O(w)$ and $O(w')$ in the proof of Lemma 6.}
\end{figure}

{\bf Lemma 6.} {\it If $O(w)$ is a octagonal chain with some integer $i\in\{1,2,\ldots,n-2\}$ such that $w_i=1$ or $3,$ then let $w'=\{w_1,\ldots,w_{i-1},2,4-w_{i+1},\ldots,4-w_{n-2}\}$
 and we have $K\!f(O(w))<K\!f(O(w')).$}

{\bf Proof.} The proof of the case $w_i=1$ is similar to the case $w_i=3$, so we only prove the case $w_i=3$. We first choose four vertices $u,v,x$ and $y$ in the $i+1$-th octagon of $O(w)$ as shown in Fig. 5. Deleting edges $\{ux\,,vy\}$ and adding edges $\{uy\,,vx\},$ we can get graph $O(w').$ Obviously, $O(w)$ and $O(w')$ are pairs of $S,T$-isomers. Let $A_2$ be the component of $O(w)-\{ux\,,vy\}$ such that $\{u,v\}\in V(A_2)$ and let $B_2$ be the component of $O(w)-\{ux\,,vy\}$ such that $\{x,y\}\in V(B_2).$ We denote by $u_1$ the unique neighbor of $u$ in $A_2$ and denote by $u_2$ the other neighbor of $u_1.$ Let $x_1$ the unique neighbor of $x$ in $B_2;$ See Fig. 5. Note that if $k\in V(A_2)\backslash \{u,u_1,v\},$ we have $r_{A_2}(u,k)=r_{A_2}(u_2,k)+2$ and $r_{A_2}(v,k)=r_{A_2}(v_1,k)+1.$ By Lemma 4, We have
$$
K\!f(O(w))-K\!f(O(w'))=\frac{[r_{A_2}(u)-r_{A_2}(v)][r_{B_2}(y)-r_{B_2}(x)]}{r_{A_2}(u,v)+r_{B_2}(x,y)+2}.
$$
We distinguish two cases. If $k\in V(A_3)\backslash \{u,u_1,u_2,v,v_1\},$ we have
\begin{equation}
  r_{A_2}(v,k)\leq 1+r_{A_2}(v_1,u_2)+r_{A_2}(u_2,k)< 2+r_{A_2}(u_2,k)< r_{A_2}(u,k).
\end{equation}

If $k\in \{u_1,u_2,v_1\},$ we have
\begin{equation}
\begin{split}
 &r_{A_2}(v,v_1)+r_{A_2}(v,u_2)+r_{A_2}(v,u_1)=4+2r_{A_1}(v,u_2)\\
 <& 5+ r_{A_1}(v,u_2)=r_{A_2}(u,v_1)+r_{A_2}(u,u_2)+r_{A_2}(u,u_1).
 \end{split}
\end{equation}

By the definitions of $r_{A_2}(u)$ and $r_{A_2}(v),$ we have
$$
r_{A_2}(u)-r_{A_2}(v)=\sum_{k\in V(A_2)\backslash\{u,v\}}r_{A_2}(u,k)-\sum_{k\in V(A_2)\backslash\{u,v\}}r_{A_2}(v,k).
$$
Thus, by equation (1)-(2), we have $r_{A_2}(u)-r_{A_2}(v)>0.$

Next, we will proof $r_{B_2}(y)-r_{B_2}(x)<0.$ Let $k\in V(B_2)\backslash \{x,y\},$ we have
$$
r_{B_2}(y,k)\leq r_{B_2}(y,x_1)+r_{B_2}(x_1,k)< 1+r_{B_2}(x_1,k)< r_{B_2}(x,k).
$$
Since
$$
r_{B_2}(y)-r_{B_2}(x)=\sum_{k\in V(B_2)\backslash\{x,y\}}r_{B_2}(y,k)-\sum_{k\in V(B_2)\backslash\{x,y\}}r_{B_2}(c,k),
$$
we deduce that $r_{B_2}(y)-r_{B_2}(x)<0.$ Thus
$$
K\!f(O(w))-K\!f(O(w'))=\frac{[r_{A_2}(u)-r_{A_2}(v)][r_{B_2}(y)-r_{B_2}(x)]}{r_{A_2}(u,v)+r_{B_2}(x,y)+2}<0.
$$
This completes the proof. \hfill $\Box$

Now we are ready to give a proof of Theorem 1.

{\bf Proof of Theorem 1.} Let $O(w)$ be an octagonal chain with the maximum Kirchhoff index. We assert that $O(w)=O(\underbrace{2,2,\ldots,2}_{n-2}).$  To the contrary, there exists some $i\in \{1,2,\ldots,n-2\}$ such that $w_i\neq 2.$ Note that $w_i=\{0,1,2,3,4\}.$ We distinguish two cases. If $w_i=0$ or $4.$  By lemma 5, we can obtain an octagonal chain with larger Kirchhoff index than $O(w),$ a contradiction.  If $w_i=1$ or $3,$  by lemma 6, we can also get a contradiction. Thus, for all $i\in \{1,2,\ldots,n-2\},$ $w_i=2.$ This completes the proof. \hfill $\Box$

{\bf Lemma 7.} {\it  If $O(w)$ is a octagonal chain with some integer $i\in\{1,2,\ldots,n-2\}$ such that $w_i=0$ or $4,$ then let $w'=\{w_1,\ldots,w_{i-1},3,4-w_{i+1},\ldots,4-w_{n-2}\}$ and
 we have $$K\!f(O(w))<K\!f(O(w')).$$}

\begin{figure}[!ht]
\centering
  \includegraphics[width=100mm]{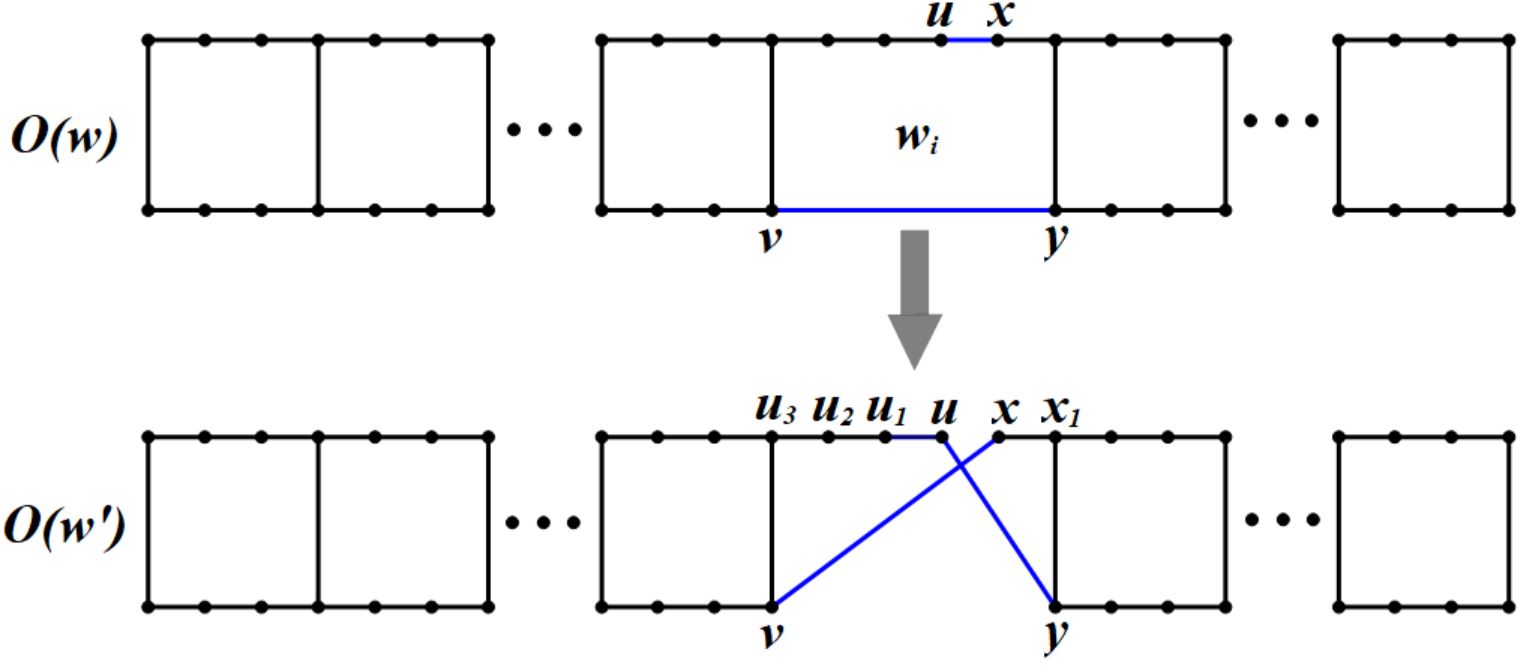}
  \caption{$O(w)$ and $O(w')$ in the proof of Lemma 7.}
\end{figure}

{\bf Proof.} The proof is very similar to the proof of Lemma 5. So we omit some details. We first prove the case $w_i=4$, and the case of $w_i=0$ is similar, so we omit the process. First choosing four vertices $u,v,x$ and $y$ in the $(i+1)$-th octagon of $O(w)$ as shown in Fig. 6. Deleting edges $\{ux\,,vy\}$ and adding edges $\{uy\,,vx\},$ we can get graph $O(w').$ Obviously, $O(w)$ and $O(w')$ are pairs of $S,T$-isomers. Let $A_3$ be the component of $O(w)-\{ux\,,vy\}$ such that $\{u,v\}\in V(A_3)$ and let $B_3$ be the component of $O(w)-\{ux\,,vy\}$ such that $\{x,y\}\in V(B_3).$ Let $u_3,u_2,u_1$ be the vertices in the $i+1$-th octagon of $O(w)$ as shown in Fig. 6. By Lemma 4, We have
$$
K\!f(O(w))-K\!f(O(w'))=\frac{[r_{A_3}(u)-r_{A_3}(v)][r_{B_3}(y)-r_{B_3}(x)]}{r_{A_3}(u,v)+r_{B_3}(x,y)+2}.
$$
Now, we are going to compare $r_{A_3}(u)$ and $r_{A_3}(v)$. Note that if $k\in V(A_3)\backslash \{u,u_1,u_2\},$ we have $r_{A_3}(u,k)=r_{A_3}(u_2,k)+2.$ We distinguish two cases. If $k\in V(A_3)\backslash \{u,u_1,u_2,u_3,v\},$ we have
$$
r_{A_3}(v,k)\leq r_{A_3}(v,u_3)+r_{A_3}(u_3,k)< 1+r_{A_3}(u_3,k)< r_{A_3}(u,k).
$$
If $k\in \{u_1,u_2,u_3\},$ we have
$$
r_{A_3}(v,u_1)+r_{A_3}(v,u_2)+r_{A_3}(v,u_3)=3+3r_{A_3}(v,u_2)< 6 =r_{A_3}(u,u_1)+r_{A_3}(u,u_2)+r_{A_3}(u,u_3).
$$
Since
$$
r_{A_3}(u)-r_{A_3}(v)=\sum_{k\in V(A_3)\backslash\{u,v\}}r_{A_3}(u,k)-\sum_{k\in V(A_3)\backslash\{u,v\}}r_{A_3}(v,k),
$$
we deduce that $r_{A_3}(u)-r_{A_3}(v)>0.$

Next, we will proof $r_{B_3}(y)-r_{B_3}(x)<0.$ Denote by $x_1$ the unique neighbor of $x$ in $B_3;$ See Fig. 6. Let $k\in V(B_3)\backslash \{x,y\},$ we have
$$
r_{B_3}(y,k)\leq r_{B_3}(y,x_1)+r_{B_3}(x_1,k)< 1+r_{B_3}(x_1,k)< r_{B_3}(x,k).
$$
Since
$$
r_{B_3}(y)-r_{B_3}(x)=\sum_{k\in V(B_3)\backslash\{x,y\}}r_{B_3}(y,k)-\sum_{k\in V(B_3)\backslash\{x,y\}}r_{B_3}(c,k),
$$
we deduce that $r_{B_3}(y)-r_{B_3}(x)<0.$ Thus
$$
K\!f(O(w))-K\!f(O(w'))=\frac{[r_{A_3}(u)-r_{A_3}(v)][r_{B_3}(y)-r_{B_3}(x)]}{r_{A_3}(u,v)+r_{B_3}(x,y)+2}<0.
$$
This completes the proof. \hfill $\Box$

For Convenience, we denote the octagons of $O(w)$ by $H_0,H_1,\ldots,H_{n-1}$ such that $H_{i-1}$ is adjacent to $H_{i}(1\leq i\leq n-1).$ Moreover, let $t_ib_i$ be the common edge of $H_{i-1}$ and $H_{i}$ such that $t_i$ is the top common vertices; See Fig. 7(i).
\begin{figure}[!ht]
\centering
  \includegraphics[width=150mm]{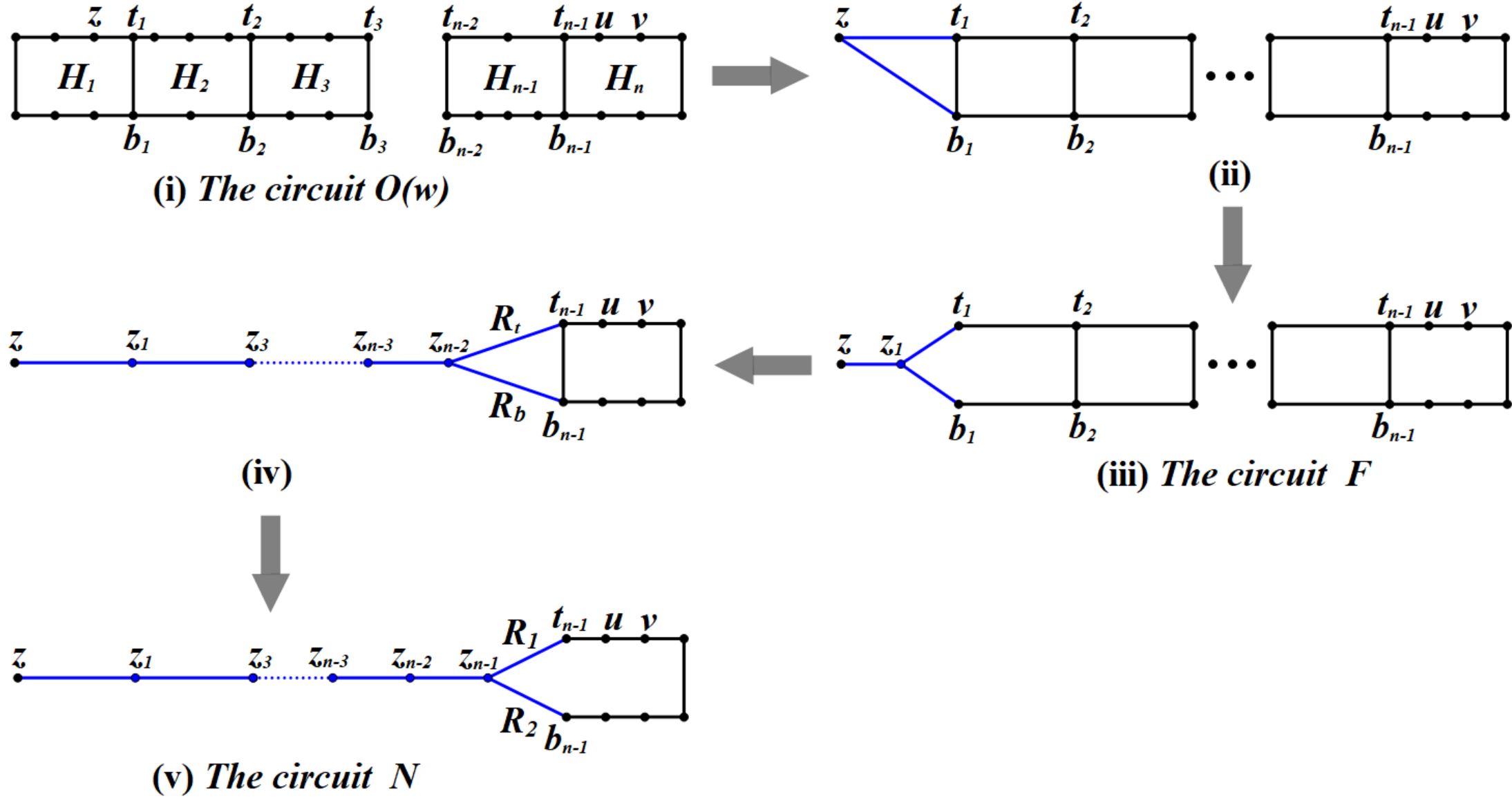}
  \caption{Illustration of circuit simplification to $O(w)$ in the proof of Lemma 8.}
\end{figure}

{\bf Lemma 8.} {\it  Let $u$ be the vertex in $H_{n-1}$ adjacent to $t_{n-1}$ with degree $2$ and let $v$ be the other neighbor of $u$. If $O(w)$ is a weighted octagonal chain and the weight on edge $t_{n-1}b_{n-1}$ is $1,$ then for any $z\in V(H_0)\backslash \{t_1,b_1\},$ we have $r_{O(w)}(z,u)<r_{O(w)}(z,v).$}

{\bf Proof.} Note that $V(H_0)\backslash \{t_1,b_1\}$ has six different vertices, we first choose $z$ the vertex which is adjacent to $t_1.$ The other remaining vertices can be proved in the same way, so we omit the process. In order to obtain our result, we need the following steps to simplify the circuit $O(w)$.

$\bullet$ First, perform the Series transform on $H_0$ of $O(w)$ to turn it into a triangle as shown in Fig. 7(ii).

$\bullet$ Next, perform the $\Delta$-$Y$ transform on this new triangle. This results in a new vertex $z_1$ as shown in Fig. 7(iii). We denote this resulting circuit $F.$

By Lemma 3, we have $$r_{O(w)}(z,u)=r_{F}(z,u),\,\,\,\, r_{O(w)}(z,v)=r_{F}(z,v).$$

By repeatedly using this two steps, we obtain the simplified circuit of $O(w)$ as depicted in Fig. 7(iv). Note that the edges $z_{n-2}t_{n-1}$ and $z_{n-2}b_{n-1}$ in Fig. $5(iii)$ are the new edges after the transformation. Denote the weighes of $z_{n-2}t_{n-1}$ and $z_{n-2}b_{n-1}$ by $R_t$ and $R_b.$ By making $\Delta$-$Y$ transformation to triangle $z_{n-2}t_{n-1}b_{n-1},$ we could obtain a simplified circuit $N,$ as shown in Fig. 7(v). Suppose that the weights of edges  $z_{n-1}t_{n-1}$ and $z_{n-1}k_{n-1}$ in $N$ are $R_1$ and $R_2;$ See Fig. 7(v). Recall that the weight of edge $t_{n-1}b_{n-1}$ is $1.$ By Definition 3, we have
\begin{eqnarray*}
  R_1=\frac{R_b}{R_a+R_b+1}, \ \  R_2=\frac{R_a}{R_a+R_b+1}.
\end{eqnarray*}
Obviously, $0<R_1<1.$ Then using the parallel and series circuit reductions yields
\begin{eqnarray*}
  &r_{O(w)}(z,u)=r_N(z,u)=r_N(z,z_{n-1})+\frac{(R_1+1)(R_2+6)}{R_1+R_2+7}, \\
  &r_{O(w)}(z,v)=r_N(z,v)=r_N(z,z_{n-1})+\frac{(R_1+2)(R_2+5)}{R_1+R_2+7}.
\end{eqnarray*}
Hence, we have
\begin{equation*}
  r_{O(w)}(z,u)-r_{O(w)}(z,v)=\frac{(R_1+1)(R_2+6)}{R_1+R_2+7}-\frac{(R_1+2)(R_2+5)}{R_1+R_2+7}=\frac{R_1-R_2-4}{R_1+R_2+7}<0.
\end{equation*}
The third inequality follows from the condition $0<R_1<1.$ Therefore,we have $$r_{O(w)}(z,u)<r_{O(w)}(z,v).$$ This completes the proof. \hfill $\Box$

Now we are ready to give a proof of Theorem 2.
\begin{figure}[!ht]
\centering
  \includegraphics[width=100mm]{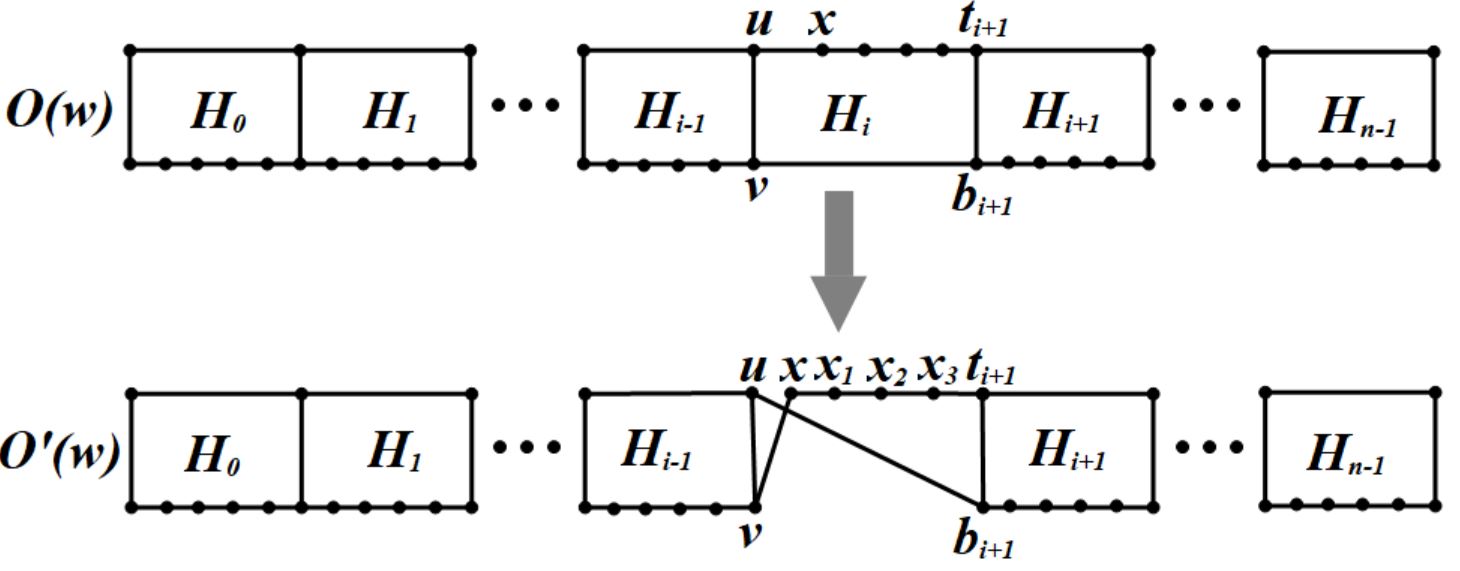}
  \caption{Illustration of circuit simplification to $O(w)$ in the proof of Theorem 2.}
\end{figure}

{\bf Proof of Theorem 2.} Let $O(w)$ be an octagonal chain containing $n$ octagons with $w=(w_1,w_2,\ldots,w_n)$. Suppose that $O(w)$ has minimum Kirchhoff index. By Lemma 5 and 7, we have $w_i=0$ or $4.$ Without loss of generality, we assume that $w_1=0.$ Now we are going to proof $w_i=0$ for all $i\in \{1,2,\ldots,n-2\}.$ If there exists some $i$ such that $w_{i-1}=0$ and $w_i=4,$ denote the vertices of the octagon $H_i$ by $u,x,x_1,x_2,x_3,t_{i+1},b_{i+1},v$ as show in Fig. 8. Let $O(w')$ be the graph obtained from $O(w)$ by deleting edges $\{ux,vb_{i+1}\}$ and adding two new edges $\{ub_{i+1},vx\}.$ Next, we will proof $K\!f(O(w'))<K\!f(O(w)).$ Then we will get a contradiction.

By the construction of $O(w'),$ we have that $O(w)$ and $O(w')$ are pairs of $S,T$-isomers. Let $A$ be the component of $O(w)-\{ux,vb_{i+1}\}$ such that $\{u,v\}\subseteq V(A),$ and  Let $B$ be the component of $O(w)-\{ux,vb_{i+1}\}$ such that $\{x,v\}\subseteq V(A).$ By Lemma 4, we have
\begin{equation}
K\!f(O(w))-K\!f(O(w'))=\frac{[r_A(u)-r_A(v)][r_B(b_{i+1})-r_B(x)]}{r_A(u,v)+r_B(x,b_{i+1})+2}.
\end{equation}

We first consider $r_A(u)$ and $r_A(v).$ Let $z$ be a vertex in $A,$ we distinguish two cases.

Case 1.  $k\in V(H_j),0\ \leq j\leq i-2.$ By using series and parallel connection rules, we can simplify $A$ to a weighted octagonal chain consisting of octagons $H_j, H_{j+1}, \ldots,H_{i-1}.$ Note that the weight on edge $t_{i-1}b_{i-1}$ is 1. by Lemma 8, we have
$r_A(k,u)<r_A(k,v).$

Case 2. $k\in V(H_{i-1}).$ By using series and parallel connection rules, we can simplify $A$ to a weighted octagon $H_{i-1}$ with the weight $r$ on edge $t_{i-1}b_{i-1}$ and the weight $1$ on other edges. Then
\begin{eqnarray*}
  &\sum_{k\in V(H_{i-1})}r_A(k,u)
  =\frac{r+6}{r+7}+\frac{6(r+1)}{r+7}+\frac{5(r+2)}{r+7}+\frac{4(r+3)}{r+7}+\frac{3(r+4)}{r+7}+
  \frac{2(r+5)}{r+7}+\frac{r+6}{r+7},\\
  &\sum_{k\in V(H_{i-1})}r_A(k,v)
  =\frac{r+6}{r+7}+\frac{2(r+5)}{r+7}+\frac{5(r+2)}{r+7}+\frac{4(r+3)}{r+7}+\frac{3(r+4)}{r+7}+
  \frac{2(r+5)}{r+7}+\frac{r+6}{r+7}.
\end{eqnarray*}
Noting that the initially weight of edge $t_{i-1}b_{i-1}$ is $1,$ we have $r<1.$ Since
\begin{equation*}
  \sum_{k\in V(H_{i-1})}r_A(k,u)-\sum_{k\in V(H_{i-1})}r_A(k,v)=\frac{4r-4}{r+7}<0,
\end{equation*}
we have
\begin{equation*}
 \sum_{k\in V(H_{i-1})}r_A(k,u)<\sum_{k\in V(H_{i-1})}r_A(k,v).
\end{equation*}
By Case 1-2, we deduce that $r_A(u)<r_A(v).$

We then consider $r_B(x)$ and $r_B(y).$ Note that if $k\in V(B)\backslash \{x,x_1,x_2,x_3\},$ we have $r_{B}(x,k)=r_{B}(t_{i+1},k)+4.$ We distinguish two cases. If $k\in V(B)\backslash \{x,x_1,x_2,x_3,t_{i+1}\},$ we have
$$
r_{B}(y,k)\leq r_{B}(y,t_{i+1})+r_{B}(t_{i+1},k)< 1+r_{B}(t_{i+1},k)< r_{B}(x,k).
$$
If $k\in \{x_1,x_2,x_3,t_{i+1}\},$ we have
\begin{align*}
&r_{B}(y,x_1)+r_{B}(y,x_2)+r_{B}(y,x_3)+r_{B}(y,t_{i+1})=6+4r_{B}(v,u_2)<10 \\
=&r_{B}(x,x_1)+r_{B}(x,x_2)+r_{B}(x,x_3)+r_{B}(x,t_{i+1}).
\end{align*}
Since
$$
r_{B}(y)-r_{B}(x)=\sum_{k\in V(B)\backslash\{u,v\}}r_{B}(y,k)
-\sum_{k\in V(B) \backslash\{u,v\}}r_{B}(x,k),
$$
we deduce that $r_{B}(y)-r_{B}(x)>0.$ By equation (1), we obtain $K\!f(O(w))>K\!f(O(w')).$
This completes the proof. \hfill $\Box$

\end{document}